\documentclass[11pt]{article}
\usepackage{epsfig,latexsym,amsfonts,amsmath,amssymb
}

\catcode`\@=11
\@addtoreset{equation}{section}

\catcode`\@=12

\def\textsf{\bf}  

\newtheorem{Theorem}{Theorem}
\newtheorem{Lemma}
{Lemma}
\newtheorem{Proposition}
{Proposition}
{Corollary}
\newtheorem{Remark}
{Remark}

\def\proof{\noindent{\textbf{Proof. }}}

\def\tildeSobolev{\widetilde{H}^m(\Omega)}

\def\Dspace{\mathcal D^m(\R^n)}

\def\tildeSobolev{\widetilde{H}^m(\Omega)}
\def\tildeSobolevN{\widetilde{H}^m_N(\Omega)}

\def\Dm{\left(-\Delta\right)_D^{\!m}\!}                                  

\def\Ds{\left(-\Delta\right)_D^{\!m}\!}
\def\Dshalf{\left(-\Delta\right)_D^{\!\frac{m}{2}}\!}

\def\Dalpha{\left(-\Delta\right)_D^{\!\alpha}\!}

\def\DalphaO{\left(-\Delta\right)_D^{\!\alpha}\!}

\def\DsON{\left(-\Delta\right)_N^{\!m}\!}                 
\def\DshalfON{\left(-\Delta\right)_N^{\!\frac{m}{2}}\!}
\def\DalphaON{\left(-\Delta\right)_{N}^{\!\alpha}\!}

\def\QED{\hfill {$\square$}\goodbreak \medskip}

\def\f{\varphi}
\def\irn{\int\limits_{\R^n}}
\def\div{{\rm div}}

\def\R{\mathbb R}

\def\pstars{{p^*_m}}

\begin{document}

\title 
{On the Sobolev and Hardy constants\\ for the fractional Navier Laplacian}

\author{Roberta Musina\footnote{Dipartimento di Matematica ed Informatica, Universit\`a di Udine,
via delle Scienze, 206, 33100 Udine, Italy. Email: {roberta.musina@uniud.it}. {Partially supported by Miur-PRIN 201274FYK7\_004,
``Variational and perturbative aspects of nonlinear differential problems''.}
}~ and
Alexander I. Nazarov\footnote{St.Petersburg Department of Steklov Institute, Fontanka, 27, St.Petersburg, 191023, Russia
and St.Petersburg State University, 
Universitetskii pr. 28, St.Petersburg, 198504, Russia. E-mail: {al.il.nazarov@gmail.com}. Supported by RFBR grant 14-01-00534 and by
St.Petersburg University grant 6.38.670.2013.}
} 

\date{}

\maketitle

\footnotesize

\noindent
{\bf Abstract.}
We prove the coincidence of the Sobolev and Hardy constants relative to 
 the ``Dirichlet''  and ``Navier'' fractional Laplacians of any real order
$m\in(0,\frac{n}{2})$ over bounded domains in $\R^n$.
\normalsize

\bigskip

\bigskip

\section{Introduction}
For any integer $n\ge 1$ the 
(fractional) Laplacian of real order $m>0$ over $\R^n$ is defined by
$$
\mathcal F\left[\Ds u\right]=|\xi|^{2m}\mathcal F[u]~\!,
$$
where $\mathcal F$ is the Fourier transform
$$
{\cal F}[u](\xi)=\frac{1}{(2\pi)^{n/2}}\irn e^{-i~\!\!\xi\cdot x}u(x)~\!dx~\!.
$$
Let $p\in(1,\infty)$ and assume  $n>pm$. Put $I_m(f)=|x|^{m-n}\star f~\!$.
Then the Hardy--Littlewood--Sobolev inequality \cite{HL1, HL2, So} 
states that $I_m$ is continuous operator
from $L^p(\R^n)$ to $L^{p^*_m}(\R^n)$, where
$$
p^*_m:=\frac{pn}{n-pm}
$$
is the {\em critical Sobolev exponent}.

We denote by $\mathcal D^{m,p}(\R^n)$ the image of $I_m$. Since for any  $f\in L^p(\R^n)$
$$
\Dshalf\big(|x|^{m-n}\star f\big)=c_{n,m}\cdot f
$$
in the distributional sense on $\R^n$ (here the constant  $c_{n,m}$ depends only on $n$ and $m$), we have 
$$
\mathcal D^{m,p}(\R^n)
=\{u\in L^\pstars(\R^n)~|~\Dshalf u\in L^p(\R^n)~\}.
$$
We endow $\mathcal D^{m,p}(\R^n)$ with the norm
$$
\|u\|_{\mathcal D^{m,p}}=\|\Dshalf u\|_p:=\Big(\int\limits_{\R^n}|\Dshalf u|^p~\!dx\Big)^{1/p},
$$
so that $I_m: L^p(\R^n)\to \mathcal D^{m,p}(\R^n)$ is
(up to a constant) an isometry with inverse $\Dshalf$~\!. In particular, $\mathcal D^{m,p}(\R^n)$ is a reflexive Banach space.

In the Hilbertian case $p=2$ we will simply write $\Dspace$ instead of $\mathcal D^{m,2}(\R^n)$. 
The explicit value and the extremals of the best constant ${\cal S}_m$ in the inequality 
$$
\irn|\Dshalf u|^2~\!dx\ge {\cal S}_m\Big(\irn|u|^{2^*_m}~\!dx\Big)^{\!\!\frac{2}{2^*_m}}
\quad\text{for any $u\in \mathcal D^{m}(\R^n)$}
$$
were furnished by Cotsiolis and Tavoularis in \cite{CoTa}. 

Next, we introduce the ``Dirichlet'' Laplacian of order $m$ over a bounded and smooth domain $\Omega\subset\R^n$ via the quadratic form
$$
Q_m^D[u]=(\Dm u,u):=\int\limits_{\mathbb R^n}|\Dshalf u|^2~\! dx~\!,
$$
with domain 
$$
\tildeSobolev=\{u\in \mathcal D^m(\R^n)\,:\,{\rm supp}\, u\subset\overline{\Omega}\}~\!.
$$
We endow $\tildeSobolev$ with the norm $\|~\!\cdot~\!\|_{\mathcal D^{m}}$. Since ${\cal C}^\infty_0$ is dense in $\mathcal D^m(\R^n)$,
a standard dilation argument implies that
$$
{\cal S}_m=\inf_{\scriptstyle u\in\widetilde H^m(\Omega)\atop u\neq 0}\frac
{Q^D_m[u]}{\|u\|_{2^*_m}^2}~\!.
$$

We introduce also the ``Navier'' Laplacian $\DsON$ of order $m$ over $\Omega$ as the $m^{\text {th}}$ power of the conventional Laplacian 
$-\Delta$ on $H^1_0(\Omega)$,
in the sense of spectral theory. More precisely, for $u\in L^2(\Omega)$ we define
$$
\DsON u :=\sum_{j\ge 1}\lambda_j^m\Big(\int\limits_\Omega u\varphi_j~\!dx\Big)\varphi_j.
$$
Here $\lambda_j, \f_j$ are, respectively, the eigenvalues and eigenfunctions (normalized in $L^2(\Omega)$) of $-\Delta$ on $H^1_0(\Omega)$
while the series converges in the sense of distributions.

The corresponding quadratic form is
$$
Q^N_m[u]=((-\Delta)^m_Nu,u)=\sum_{j\ge 1}\lambda_j^m\Big(\int\limits_\Omega u\varphi_j~\!dx\Big)^2
=\int\limits_\Omega|\DshalfON u|^2~\!dx,
$$
with domain 
$$
\tildeSobolevN=\{u\in L^2(\Omega)\,:\, Q^N_m[u]<\infty\}.
$$ 

Finally, we
define the {\em Navier-Sobolev} constant by
$$
{\cal S}_m^N:=\inf_{\scriptstyle u\in\widetilde H^m_N(\Omega)\atop u\neq 0}\frac
{Q^N_m[u]}{\|u\|_{2^*_m}^2}~\!.
$$
We are in position to state the main result of the present paper.
\begin{Theorem}
\label{T:Sobolev_constants}
Let $\Omega$ be a bounded and smooth domain in $\R^n$ and $m\in\big(0,\frac{n}{2}\big)$. Then
$$
{\cal S}_m^N= {\cal S}_m.
$$
\end{Theorem}

Our argument applies also to Hardy-Rellich type inequalities. 
The explicit value of the positive constant 
$$
{\cal H}_m:=\inf_{\scriptstyle u\in\Dspace\atop U\neq 0}\frac{Q^D_m[u]}{\||x|^{-m}u\|_{2}^2}
=\inf_{\scriptstyle u\in\tildeSobolev\atop U\neq 0}\frac{Q^D_m[u]}{\||x|^{-m}u\|_{2}^2}
$$
has been computed in \cite{He} (see also \cite{DH} and \cite{Mit2} for the integer orders $m\in\mathbb N$, even 
in a non-Hilbertian setting). The {\em Navier-Hardy} constant over 
a bounded and smooth domain $\Omega$ is defined by
$$
{\cal H}_m^N:=\inf_{\scriptstyle u\in\widetilde H^m_N(\Omega)\atop u\neq 0}\frac
{Q^N_m[u]}{\||x|^{-m}u\|^2_{2}}~\!.
$$
The argument we use to prove Theorem \ref{T:Sobolev_constants} plainly
leads to the next result.
\begin{Theorem}
\label{T:Hardy_constants} 
Let $\Omega$ be a bounded and smooth domain in $\R^n$ and $m\in\big(0,\frac{n}{2}\big)$. Then
$$
{\cal H}_m^N={\cal H}_m.
$$
\end{Theorem}

The equalities ${\cal S}^N_1={\cal S}_1$,  ${\cal H}^N_1={\cal H}_1$ are totally trivial. 
If $m\neq 1$ is an integer number, then the inequalities ${\cal S}^N_m\le {\cal S}_m$ and ${\cal H}^N_m\le {\cal H}_m$ follow immediately from 
$\tildeSobolev\subseteq\widetilde H^m_N(\Omega)$, whereas the opposite inequalities need a detailed proof.

For integer orders $m\in \mathbb N$, the statements of Theorems \ref{T:Sobolev_constants}
and \ref{T:Hardy_constants} are known (even in non-Hilbertian setting).
The coincidence of the two Hardy constants
can be extracted from the proof of Theorem 3.3 in \cite{Mit2} (see also \cite[Lemma 1]{GGM}),
where Enzo Mitidieri took advantage of a Rellich--Pokhozhaev type identity \cite{P, MitId}. 
The coincidence of the two Sobolev  constants for  $m\in\mathbb N$  was obtained in \cite{GGS}
(see also \cite{Ge, VdV} for previous results in case $p=2$ and $m=2$). 
We cite also \cite{M2}, 
where weighted Sobolev constants are studied under the hypothesis $m=2$. 

We emphasize that for $m\notin\mathbb N$ none of the inequalities
${\cal S}^N_m\le {\cal S}_m$, ${\cal S}^N_m\ge {\cal S}_m$ (respectively, ${\cal H}^N_m\le {\cal H}_m$, ${\cal H}^N_m\ge {\cal H}_m$) 
is easily checked. For $m\in(0,1)$, Theorem \ref{T:Sobolev_constants} was proved in 
\cite{FL}. To handle the general case of real orders $m>0$ we largely use some of the results in \cite{FL, FL2}. 
Additional tools are the maximum principles for fractional Laplacians
and a result about the transform $u\mapsto |u|$, $u\in \tildeSobolev$, for $0<m<1$, that
might have an independent interest (see Theorem \ref{T:truncation}).

\section{Preliminaries}
Here we collect some facts about the Dirichlet and the Navier quadratic forms.
\paragraph{1.} First, we note that $\tildeSobolev\subseteq\widetilde H^m_N(\Omega)$ and
$$
\tildeSobolev=\widetilde H^m_N(\Omega)\quad\text{if and only if~ $m<\displaystyle\frac{3}{2}$}.
$$
This fact is well known for natural orders $m$; the general case  follows immediately
from \cite[Theorem 1.17.1/1]{Tr} and \cite[Theorem 4.3.2/1]{Tr}. 
\paragraph{2.} It is well known that for any $m\in\mathbb N$
$$
\widetilde H^m_N(\Omega)=
\displaystyle{\left\{u\in H^m(\Omega)~\Big|~
\text{tr}_{\partial\Omega}\big[(-\Delta)^\nu u\big]=0\ \ 
\text{for}\ \ \nu\in\mathbb N_0, ~ \nu<\frac{m}{2}
\right\}}.
$$
We omit the proof of the next simple analog for non integer $m$.

\begin{Lemma}
\label{L:spaces}
Let $m\notin\mathbb N$, $m>1$.
\begin{itemize}
 \item 
 If $\lfloor m\rfloor\ge 2$ is even, then $\widetilde H^m_N(\Omega)=\Big\{u\in \widetilde H^{\lfloor m\rfloor}_N(\Omega)~\Big|
~(-\Delta)_N^{\!\frac{\lfloor m\rfloor}{2}}u\in \widetilde H^{m-\lfloor m\rfloor\!}(\Omega)~\Big\}$.

 \item 
If $\lfloor m\rfloor\ge 1$ is odd, then \ $\widetilde H^m_N(\Omega)=\Big\{u\in \widetilde H^{\lfloor m\rfloor}_N(\Omega)~\Big|
~(-\Delta)_N^{\!\frac{m}{2}}u\in L^2(\Omega)~\Big\}$.

\end{itemize}
\end{Lemma}

\paragraph{3.} Let $m\in\mathbb N$ and let $u\in \tildeSobolev$. Then it is easy to see that 
$Q^D_m[u]= Q^N_m[u]$. More precisely, if $m$ is even one gets the pointwise equality
$$
\Dshalf u=\DshalfON u=(-\Delta)^{\!\frac{m}{2}} u.
$$
If $m$ is odd the following integral equalities hold:
$$
\irn|\Dshalf u|^2~\!dx=\int\limits_\Omega|\DshalfON u|^2~\!dx=\int\limits_\Omega|\nabla(\Delta^{\!\frac{m-1}{2}} u)|^2~\!dx.
$$
Integrating by parts we can write 
for all $m\in\mathbb N$
\begin{equation}
\label{gradient}
\irn|\Dshalf u|^2~\!dx=\int\limits_\Omega|\DshalfON u|^2~\!dx=\int\limits_\Omega|\nabla^m u|^2~\!dx,
\qquad u\in \tildeSobolev.
\end{equation}
For non integer orders $m$ the Dirichlet and Navier quadratic forms never
coincide on the Dirichlet domain $\tildeSobolev$. Indeed, the next result holds.

\begin{Proposition}[\cite{FL, FL2}]
\label{T:forms}
Let $m>0$, $m\notin\mathbb N$, and let $u\in\tildeSobolev$, $u\not\equiv 0$. 
Then
\begin{eqnarray*}
\irn|\Dshalf u|^2&<&\int\limits_\Omega|\DshalfON u|^2~\!dx \quad\text{\rm if \ $\lfloor m\rfloor$ is even;}\\
\irn|\Dshalf u|^2&>&\int\limits_\Omega|\DshalfON u|^2~\!dx\quad\text{\rm if \ $\lfloor m\rfloor$ is odd.}
\end{eqnarray*}
\end{Proposition}

In view of Proposition \ref{T:forms}, one is lead to ask ``how much'' the Dirichlet and Navier quadratic forms
differ on $\tildeSobolev$ if $m\notin\mathbb N$. The answer  takes into account  the action of dilations.

Fix any point $x_0\in\Omega$ and take $u\in\tildeSobolev$. Concentrate $u$ around $x_0$ by putting
$u_\rho(x)=\rho^{\frac{n-2m}{2}}u(\rho(x-x_0)+x_0)$ for $\rho\gg 1$. Then
$u_\rho\in\tildeSobolev$ and  $Q^D_m[u_\rho]\equiv Q^D_m[u]$. In contrast,
$Q^N_m[u_\rho]$ depends on $\rho$, as the Navier quadratic form does depend on the domain $\Omega$. 
Nevertheless, the next result holds.

\begin{Proposition}[\cite{FL, FL2}]
\label{T:3}
Let $m>0$ and  $u\in\tildeSobolev$. 
Then 
$$
\irn|\Dshalf u|^2~\!dx=\lim_{\rho\to \infty}\int\limits_\Omega|\DshalfON u_\rho|^2~\!dx.
$$
\end{Proposition}

\paragraph{4.} It is well known that if $u\in \widetilde H^1(\Omega)=\widetilde H^1_N(\Omega)=H^1_0(\Omega)$
then $|u|\in \widetilde H^1(\Omega)$, and $|\nabla |u||=|\nabla u|$ almost everywhere on $\Omega$. 
By (\ref{gradient}), this implies
$$\int\limits_{\R^n}|\left(-\Delta\right)_D^{\!\frac{1}{2}}\! |u||^2~\!dx =
\int\limits_{\R^n}|\left(-\Delta\right)_D^{\!\frac{1}{2}}\! u|^2~\!dx
=\int\limits_{\Omega}|\left(-\Delta\right)_N^{\!\frac{1}{2}}\! |u||^2~\!dx 
= \int\limits_{\Omega}|\left(-\Delta\right)_N^{\!\frac{1}{2}}\! u|^2~\!dx~\!.
$$
For smaller orders $m\in(0,1)$ one still has that 
$\tildeSobolev=\tildeSobolevN$ (see point 1 above), but the 
operator $u\mapsto |u|$ behaves quite differently.

\begin{Theorem}
\label{T:truncation}
Let $m\in(0,1)$ and $u\in \tildeSobolev$. Then $|u|\in\tildeSobolev$ and
\begin{eqnarray}
\label{eq:large}
\int\limits_{\R^n}|\Dshalf |u||^2~\!dx &\le& \int\limits_{\R^n}|\Dshalf u|^2~\!dx\\
\label{eq:largeN}
\int\limits_{\Omega}|\DshalfON |u||^2~\!dx &\le& \int\limits_{\Omega}|\DshalfON u|^2~\!dx~\!.
\end{eqnarray}
In addition, if both the positive and the negative parts of $u$ are nontrivial, then  strict
inequalities hold in (\ref{eq:large}) and in (\ref{eq:largeN}).
\end{Theorem}

\proof
In the paper \cite{CaSi}, the Dirichlet fractional Laplacian of order $m\in(0,1)$ was connected with the so-called 
{\it harmonic extension in $n+2-2m$ dimensions} (see also \cite{CT} for the case $m=\frac 12$). Namely, it was shown that for any 
$v\in\tildeSobolev$, the function $w_{v}(x,y)$ minimizing the weighted Dirichlet integral
$$
{\cal E}_m(w)=\int\limits_0^\infty\!\int\limits_{\mathbb{R}^n} y^{1-2m}|\nabla w(x,y)|^2\,dxdy
$$
over the set
$${\cal W}(v)=\Big\{w(x,y)\,:\,
{\cal E}_m(w)<\infty~,\ \ w\big|_{y=0}=v\Big\},
$$
satisfies
\begin{equation}
\int\limits_{\R^n}|\DshalfON v|^2~\!dx =  {c_m}~\! {\cal E}_m(w_{v}),
\label{quad_D}
\end{equation}
where the constant  $c_m$ depends only on $m$.

For any fixed $u\in \tildeSobolev$ find $w_u\in {\cal W}(u)$ and $w_{|u|}\in {\cal W}(|u|)$.
Then clearly $|w_u|\in {\cal W}(|u|)$ and therefore 
${\cal E}_m(w_{|u|})\le {\cal E}_m(|w_u|)={\cal E}_m(w_u)$. Thus (\ref{eq:large}) holds, thanks to (\ref{quad_D}).

Now assume that $u$ changes sign. The function
$w_{|u|}(x,y)$ is the unique solution of the boundary value problem
\begin{equation}
\label{eq:CS}
-\div (y^{1-2m}\nabla w)=0\quad \mbox{in}\quad \mathbb R^n\times\mathbb R_+;\qquad w\big|_{y=0}=|u|
\end{equation}
with finite energy. Hence $w_{|u|}$ is analytic in $\R^n\times\mathbb R_+$. 
Since $w_u$ changes sign then $|w_u|$ can not solve (\ref{eq:CS}), that implies
${\cal E}_m(|w_u|)> {\cal E}_m(w_{|u|})$. Hence the strict inequality holds in (\ref{eq:large}), that
concludes the proof for the Dirichlet Laplacian.

To check (\ref{eq:largeN}) one has to use, 
instead of \cite{CaSi}, the characterization of the Navier fractional Laplacian given (among some other fractional operators) in \cite{ST}.
Namely, for any $v\in\tildeSobolev$, the function $w^N_{v}(x,y)$ minimizing 
${\cal E}_m(w)$ over the set
$${\cal W}^N(v)=\Big\{w\in {\cal W}(v)\,:\,
\text{supp}~\!w(~\!\cdot~\!,y) \subseteq \overline\Omega\quad\text{for any $y>0$}~\Big\},
$$
satisfies
$$
\int\limits_{\Omega}|\DshalfON v|^2~\!dx =  {c_m}~\! {\cal E}_m(w_{v}).
$$
The rest of the proof runs as in the Dirichlet case.
We omit details.
\QED

\begin{Remark}
\label{R:MP}
Here we deal with  maximum principles for the operators $\Ds$
and $\DsON$~\!, $m\in (0,1)$.

Let $u\in\tildeSobolev$, and let $f=\Ds u\in (\tildeSobolev)'$ be a nonnegative and nontrivial distribution. Then it is well known that
$u\ge0$ in $\Omega$. This is actually a simple corollary to Theorem \ref{T:truncation}. 
The function $u$ is characterized variationally as the unique minimizer of the energy functional
$$
J(v)=\int\limits_{\R^n}|\Dshalf v|^2~\!dx -2\,\big\langle f,v\big\rangle
$$
on $\tildeSobolev$. We have $J(|u|)\le J(u)$ by Theorem \ref{T:truncation}. This implies $u=|u|\ge 0$, as desired, by the uniqueness of the minimizer.

By the same reason, if $u\in\tildeSobolev$ and $\DsON u =f\ge0 $ then $u\ge 0$ in $\Omega$.
\end{Remark}

\paragraph{5.} 
We conclude this preliminary section by recalling a well known fact already mentioned in the Introduction. 

\begin{Proposition}
\label{P:ex}
Let $p>1$, $m>0$, $n>2mp$. Then for any $f\in L^p(\R^n)$, problem
$$
\Ds U=f; \qquad U\in\mathcal D^{2m,p}(\R^n)
$$
has a unique solution. If in addition $f\neq 0$ is nonnegative, then 
$U>0$ in $\R^n$.
\end{Proposition}

\proof 
Up to a multiplicative constant, the unique solution $U$ is explicitly given by 
$|x|^{2m-n}\star f$. The statement readily follows.
\QED

\section{Proof of Theorems \ref{T:Sobolev_constants} and \ref{T:Hardy_constants}}
Since $\tildeSobolev\subseteq\tildeSobolevN$, then clearly
$$
{\cal S}_m^N=\inf_{\scriptstyle u\in\widetilde H^m_N(\Omega)\atop u\neq 0}\frac
{Q^N_m[u]}{\|u\|_{2^*_m}^2}\le \inf_{\scriptstyle u\in\widetilde H^m(\Omega)\atop u\neq 0}\frac
{Q^N_m[u]}{\|u\|_{2^*_m}^2}~\!.
$$
Hence, ${\cal S}_m^N\le {\cal S}_m$ by Proposition \ref{T:forms}, if $2k-1\le m\le 2k$, $k\in\mathbb N$,
and by Proposition \ref{T:3}, otherwise. By the same reason, ${\cal H}_m^N\le{\cal H}_m$.
Thus, it suffices to prove the opposite inequalities 
${\cal S}_m^N\ge {\cal S}_m$ and ${\cal H}_m^N\ge{\cal H}_m$.

\medskip

Fix any nontrivial $u\in \widetilde H^m_N(\Omega)$ and extend it by the null function. 
To conclude the proof, it is sufficient to construct a function
$U\in\Dspace$ such that
\begin{eqnarray}
U &\ge& |u|\quad\text{a.e. in \ $\R^n$};\label{eq:key_a}\\
\irn|\Dshalf U|^2~\!dx &\le& \int\limits_\Omega|\DshalfON u|^2~\!dx.\label{eq:key_b}
\end{eqnarray}
We have to distinguish between two cases.

\paragraph{1.} \underline{Case $2k+1<m\le 2k+2$, for some $k\in \mathbb N_0$}.

\medskip

We use Proposition \ref{P:ex} to fix the unique positive solution $U$ of
$$
\Dshalf U=\chi_\Omega|\DshalfON u|;\qquad U\in\Dspace,
$$
where $\chi_\Omega |\DshalfON u|$ denotes the null extension
of the function $|\DshalfON u|\in L^2(\Omega)$. 
Since (\ref{eq:key_b}) trivially holds, we only have to check (\ref{eq:key_a}), that is the 
trickiest step in the whole proof.

It is  convenient to write 
$$
\frac{m}{2}=k+\alpha~,\quad\frac{1}{2}<\alpha\le 1.
$$
Since $u\in \widetilde H^m_N(\Omega)$, then for any integer $\nu=0,\cdots, k$ the function
$u_\nu:=(-\Delta)^\nu u$ belongs to $H^1_0(\Omega)$, compare with 
Lemma \ref{L:spaces}. In addition we know that $u_k\in \widetilde H^{2\alpha}_N(\Omega)$, that 
implies 
$$
u_k\in H^1_0(\Omega)~,\quad \DalphaON u_k\in L^2(\Omega)~\!.
$$
We introduce the  solutions $\widetilde w$, $w$ to
\begin{eqnarray*}
\DalphaON \widetilde w=|\DalphaON u_k|;\qquad
\widetilde w\in\widetilde H^{\alpha}(\Omega);\\
\DalphaO {w}=|\DalphaON u_k|;\qquad
{w}\in \widetilde H^{\alpha}(\Omega)~\!.
\end{eqnarray*}
We claim that
\begin{equation}
\label{eq:basic}
w\ge \widetilde w\ge |u_k|\quad\text{a.e. in \ $\Omega$.}
\end{equation}
The fact that $\widetilde w\ge |u_k|$ readily follows from the maximum principle, 
see Remark \ref{R:MP} or \cite[Lemma 2.5]{CDDS}. Also by the maximum principle
$w$ is nonnegative, and hence by \cite[Theorem 1]{FL} we have 
$\DalphaON w\ge \DalphaO {w}$ in the distributional sense on $\Omega$.
Therefore, 
$$
\DalphaON (w-\widetilde w)\ge \DalphaO w-\DalphaON\widetilde w=0,
$$ 
and the maximum principle applies again to get (\ref{eq:basic}).

\bigskip

Now we decompose $U\in \Dspace$ in the same way as we did for $u$. Namely, we define
$U_\nu=(-\Delta)^\nu U$ for any integer $\nu=0,\cdots,k$, and notice that
$$
(-\Delta)^{\frac{m}{2}-\nu}_DU_\nu= \chi_\Omega|\DshalfON u|~,\qquad
U_\nu\in \mathcal D^{m-2\nu}(\R^n).
$$
By  Proposition \ref{P:ex}, $U_\nu>0$ on $\R^n$.
In particular, the function $U_k\in \mathcal D^{2\alpha}(\R^n)$ solves
$$
(-\Delta)^\alpha_D U_k=\DalphaO w\quad\text{in} \ \ \Omega;\qquad U_k>0=w\quad\text{in} \ \R^n\setminus\overline\Omega.
$$
Therefore $U_k\ge w$ on $\Omega$, and we have by (\ref{eq:basic})
\begin{equation}
\label{eq:Uh}
U_k\ge |u_k|\quad\text{a.e. in \ $\Omega$}.
\end{equation}
If
$k=0$ then we are done. If 
$k\ge 1$ then (\ref{eq:Uh}) is equivalent to
$$
-\Delta U_{k-1}\ge |-\Delta u_{k-1}|\quad\text{a.e. in \ $\Omega$,}
$$
that readily implies $U_{k-1}\ge |u_{k-1}|$ on $\Omega$, as $U_{k-1}>0$ on $\R^n$
and $u_{k-1}\equiv 0$ on $\R^n\setminus\overline\Omega$. Repeating the same argument we arrive at 
(\ref{eq:key_a}), and the proof is complete.

\paragraph{2.}
\underline{Case $2k<m\le 2k+1$, for some $k\in\mathbb N_0$. }\\
Now we write
$$
\frac{m}{2}=k+\alpha~,\quad 0<\alpha\le \frac{1}{2}.
$$
From $u\in  \widetilde H^m_N(\Omega)$ we infer that 
$(-\Delta)^k u\in \widetilde H^{2\alpha}(\Omega)$ by Lemma \ref{L:spaces}. 
Since $2\alpha\in(0,1]$, then also $|(-\Delta)^k u|\in \widetilde H^{2\alpha}(\Omega)$. 
By Sobolev embedding, $|(-\Delta)^ku|\in L^{2^*_{2\alpha}}(\Omega)$.

Notice that 
$n>2k\cdot2^*_{2\alpha}$. Therefore we can apply Proposition \ref{P:ex} with $m=k$ and $p=2^*_{2\alpha}$ to find the 
unique positive solution $U$ to 
$$
(-\Delta)^k U=|(-\Delta)^k u|;\qquad U\in \mathcal D^{2k,2^*_{2\alpha}}(\R^n)~\!.
$$
Since $(2^*_{2\alpha})^*_{2k}=2^*_m$, the Sobolev embedding theorem gives $U\in L^{2^*_m}(\R^n)$. Moreover, from 
$(-\Delta)^k U\in \mathcal D^{2\alpha}(\R^n)$
we infer that $\Dshalf U \in L^2(\R^n)$, that is, $U\in \Dspace$.

The proof of (\ref{eq:key_a}) runs now in the same way as in the case 1, and is even more simple
since we only have to handle Laplacians of integer orders.

To check (\ref{eq:key_b}), we write
$$\aligned
\irn|\Dshalf U|^2~\!dx&=\int\limits_{\R^n}|\Dalpha\big(|(-\Delta)^k u|\big)|^2~\!dx \\
&\le\int\limits_{\R^n}|\DalphaO((-\Delta)^k u)|^2~\!dx\\
&\le\int\limits_{\Omega}|\DalphaON((-\Delta)^k u)|^2~\!dx
=\int\limits_{\Omega}|\DshalfON u)|^2~\!dx~\!.
\endaligned
$$
Here the first inequality holds by Theorem \ref{T:truncation}, the second one 
follows from (\ref{gradient}) for $2\alpha=1$ and from Proposition \ref{T:forms} for $2\alpha\in(0,1)$.

Thus, Theorems \ref{T:Sobolev_constants} and \ref{T:Hardy_constants} are completely proved.
\QED

\begin{Remark}[Non-Hilbertian case]
Let $m\in\mathbb N$, and let  $1<p<\frac{n}{m}$.
With minor modifications, one gets an alternative proof of
\cite[Theorems 1 and 2]{GGS} concerning the Navier-Sobolev and Navier-Hardy constants
for the space $W^{m,p}_N(\Omega)$. 
Best constants in weighted Sobolev inequalities
can be included as well, see \cite{M2} for $m=2$.
\end{Remark}

\end{document}